\documentclass[letterpaper, paper,10pt]{AAS}		

\PaperNumber{20-647}

\usepackage{hyperref}
\hypersetup{colorlinks=true, pdfstartview=FitV, linkcolor=black, citecolor= black, urlcolor= black}
\usepackage{overcite}
\usepackage{footnpag}			      	

\usepackage{bm}
\usepackage{amsmath}
\usepackage[per-mode = symbol]{siunitx}  
\RequirePackage{color}
\usepackage{xcolor}
\usepackage{soul}

\usepackage[section]{placeins}  

\usepackage[super]{nth}
\usepackage{mathtools}

\allowdisplaybreaks

\usepackage{caption}
\usepackage{subcaption}
\usepackage{amsthm}

\usepackage{cleveref}
\newtheorem{proposition}{Proposition}

\theoremstyle{remark}
\newtheorem{assumption}{Assumption}
\newtheorem{remark}{Remark}

\newcommand{\norm}[1]{\left\lVert#1\right\rVert}

\def \VerticalAscent {1}
\def \PitchOver {2}
\def \ZlgtOne {3}
\def \CoastingOne {4}
\def \ZlgtTwo {5}
\def \CoastingTwo {6}
\def \StageThreeOne {7}
\def \StageThreeTwo {8}
\def \CoastingThree {9}
\def \StageFourOne {10}
\def \CoastingFour {11}
\def \StageFourTwo {12}
\def \Return {13}

\def \ti {t_0}
\def \tf {t_f}

\begin{document}

\title{Convex Optimization of Launch Vehicle \protect\\ Ascent Trajectory with Heat-Flux and \protect\\ Splash-Down Constraints}

\author{Boris Benedikter\thanks{PhD Student, Department of Mechanical and Aerospace Engineering, Sapienza University of Rome, Via Eudossiana 18, Rome, 00184, Italy; boris.benedikter@uniroma1.it},
Alessandro Zavoli\thanks{Research Assistant, Department of Mechanical and Aerospace Engineering, Sapienza University of Rome, Via Eudossiana 18, Rome, 00184, Italy; alessandro.zavoli@uniroma1.it},
Guido Colasurdo\thanks{Full Professor, Department of Mechanical and Aerospace Engineering, Sapienza University of Rome, Via Eudossiana 18, Rome, 00184, Italy; guido.colasurdo@uniroma1.it}, \\
Simone Pizzurro\thanks{Research Fellow, Launchers and Space Transportation Department, Italian Space Agency, Via del Politecnico snc, Rome, 00133, Italy; simone.pizzurro@est.asi.it},
\ and Enrico Cavallini\thanks{Head of Scientific Research Department, Italian Space Agency, Via del Politecnico snc, Rome, 00133, Italy; enrico.cavallini@asi.it}
}

\maketitle{}

\begin{abstract}
This paper presents a convex programming approach to the optimization of a multistage launch vehicle ascent trajectory, from the liftoff to the payload injection into the target orbit, taking into account multiple nonconvex constraints, such as the maximum heat flux after fairing jettisoning and the splash-down of the burned-out stages.
Lossless and successive convexification are employed to convert the problem into a sequence of convex subproblems.
Virtual controls and buffer zones are included to ensure the recursive feasibility of the process and a state-of-the-art method for updating the reference solution is implemented to filter out 
undesired phenomena that may hinder convergence.
A $hp$ pseudospectral discretization scheme is used to accurately capture the complex ascent and return dynamics with a limited computational effort.
The convergence properties, computational efficiency, and robustness of the algorithm are discussed on the basis of numerical results.
{The ascent of the VEGA launch vehicle toward a polar orbit is used as case study to discuss} the interaction between the heat flux and splash-down constraints. 
Finally, a sensitivity analysis of the launch vehicle carrying capacity to different splash-down locations is presented.

\end{abstract}

\section{Introduction}

In the present state of the art,
the only propulsion system capable of providing the high thrust required for access to space is the chemical one.
However, this system allows for injecting into orbit only a small fraction of the rocket initial mass. 
Therefore, the ascent trajectory optimization process is of primary interest in order to increase the launcher capacity and reduce the overall mission cost.
Besides, the need for a reliable and efficient optimization tool is apparent in the preliminary design phases of a launch vehicle, to evaluate the performance of various configuration concepts, in the advanced pre-flight analysis, to assess the feasibility of specific mission scenarios,
and in the definition of optimization-based real-time guidance algorithms,
where computational speed and robustness are primary requirements.

The design of a rocket ascent trajectory is a complex optimal control problem (OCP), 
greatly sensitive to the optimization variables and characterized by highly nonlinear dynamics and numerous mission requirements. 
Over the years, various optimization methods have been proposed to solve the ascent problem.
Jurovics\cite{jurovics1961optimum} was one of the first to propose an indirect approach.
An indirect procedure also underlies the well-known DUKSUP optimization software\cite{spurlock2014duksup}, 
which has been extensively employed by NASA in the design of the Atlas, Titan, and Space Shuttle launch systems.
More recent work based on the indirect method includes both pure trajectory design applications\cite{colasurdo1995optimization, martinon2009numerical, casalino2014optimization} and closed-loop guidance algorithms\cite{calise1998design,lu2003closed, bonalli2020optimal}.
However, when dealing with real-world launch vehicle missions, the indirect method may be unappealing as it requires the derivation of the optimality conditions, which can be a burdensome task, and the solution of the resulting boundary value problem requires a meticulous initialization process to achieve convergence.
In addition, if path constraints are included in the formulation, an \textit{a priori} knowledge of the structure of the constrained arcs is necessary, which, in general, is hard to guess. 
Thus, direct methods are typically favored.

A wide spectrum of direct optimization software tools has been developed for solving the ascent problem,
such as
POST\cite{brauer1977capabilities}, OTIS\cite{vlases1990otis}, and ASTOS\cite{wiegand2010astos}.
Direct methods have been used even for solving problems very similar to the one addressed in this paper.
Spangelo and Well\cite{spangelo1994rocket} described a direct formulation of the ascent problem that takes into account the maximum heat flux and constrains the return of a spent stage by bounding the perigee of its osculating orbit at burnout.
Instead, later work by Weigel and Well\cite{weigel2000dual} includes a complete simulation of the return of the burned-out stage as an additional phase and, then, solves the resulting OCP via direct multiple shooting.
However, these approaches essentially consist in transcribing the continuous-time OCP into a general nonlinear programming (NLP) problem and, despite being easy to set up, this frequently leads to a solution that depends on provided first guess.
This is a burdensome drawback for the problem at hand because designing an accurate initialization may be nontrivial due to the aforementioned sensitivity of the problem.
Moreover, solving a general NLP problem is a computationally expensive task, with no guarantee on the optimality of the attained solution.

Convex optimization techniques are becoming increasingly popular for solving optimal control problems in the aerospace community\cite{liu2017survey}.
Convex optimization is a special class of mathematical programming that allows for the use of polynomial-time algorithms that provide a theoretically guaranteed optimal solution with a limited computational effort.
However, since most aerospace problems are not naturally convex, several \emph{convexification} techniques have been developed to convert a nonconvex problem into a convex one.
These methods are grouped into \emph{lossless} and \emph{successive} convexification techniques.
The former consist in exploiting either a convenient change of variables
or a suitable constraint relaxation to reformulate the problem as convex.
For example, A{\c{c}}{\i}kme{\c{s}}e and Blackmore\cite{accikmecse2011lossless}
proved that problems with a certain class of nonconvex control constraints can be posed equivalently as relaxed convex problems.
Instead, successive convexification offers a way to handle the nonconvexities that cannot be handled by lossless convexification
through linearization 
around a reference solution that is recursively updated.
Differently from lossless convexification, the successive linearization generates a sequence of approximated subproblems. 
The theoretical proof that also successive convexification leads to a (locally) optimal solution of the originally intended problem is available only under appropriate assumptions\cite{liu2014solving,mao2016successive,bonalli2019trajectory}.
Nevertheless, current research offers wide numerical evidence of the effectiveness of successive convexification over a broad spectrum of applications, including 
spacecraft rendezvous\cite{benedikter2019convexrendezvous},
proximity operations\cite{lu2013autonomous},
formation flying\cite{morgan2014model},
low-thrust transfers\cite{wang2018optimization},
rocket powered landing\cite{szmuk2016successive,sagliano2019generalized,liu2018fuel},
and atmospheric entry \cite{wang2016constrained}.
Convex optimization has been proposed also for solving the launch vehicle ascent trajectory problem.
However, successful applications are limited to simplified scenarios,
where a flat Earth is assumed\cite{zhang2019rapid}, 
atmospheric forces are neglected\cite{li2019optimal, li2020online},
or only the upper stage trajectory is optimized\cite{liu2014solving, cheng2017efficient}.

In this paper, a realistic dynamical model, which accounts for a Keplerian gravitational model and nonlinear aerodynamic forces, is considered, and the complete ascent problem is solved via convex optimization, building up on the authors' previous work\cite{benedikter2019convexascent}.
The VEGA launch vehicle is taken as case study, but the method can be easily extended to any other rocket.
VEGA is a four-stage launcher made up of three solid rocket motors and a small liquid rocket engine that performs the final orbit insertion maneuver~\cite{vega2006manual}.
The rocket configuration is such that the third stage burnout velocity is close to the orbital one, hence it ends up falling far away from the launch site.
In this respect, 
the impact point of the third stage must be predicted and actively constrained to a safe location.
The return of the other stages also requires a careful design, 
but it would draw in further safety-related requirements, specific of either the launch vehicle or the launch base, that are out of the scope of this paper, and it is thus neglected.
The accurate prediction of the splash-down location requires the inclusion of an additional return phase in the formulation that must be cast as a free-time phase to efficiently satisfy the zero-altitude terminal constraint.
Moreover, 
since it features a high-velocity object falling into the atmosphere, 
the accurate discretization of its dynamics 
may significantly increase 
the overall computational burden of the optimization.
Thus, a proper discretization scheme must be adopted to limit the problem dimension.
In this respect, 
a $hp$ discretization, based on Radau pseudospectral method\cite{garg2010unified}, is employed to obtain accurate solutions with a limited computational cost\cite{darby2011hp}.
Finally, a constraint related to the maximum heat flux that the payload can undergo once the fairing is jettisoned is included in the formulation. 
This is another nonconvex constraint to be tackled in the convexification procedure and it presents the further technical difficulty of being coupled with the splash-down constraint since moving the impact point affects the whole trajectory profile and, in particular, the encountered heat flux conditions.

The present paper employs a combination of lossless and successive convexification techniques to convert the launch vehicle ascent problem into a sequence of convex subproblems.
To enhance convergence, the update of the reference solution is devised as a weighted sum of the previously-found solutions, rather than being based only on the last one.
This approach successfully filters out oscillations in the search space and other common undesired phenomena due to the successive linearization, such as artificial unboundedness, and is thus referred to as \emph{filtering}\cite{benedikter2019convexrendezvous}.
Also, virtual controls and buffer zones are included in the formulation to prevent artificial infeasibility, and a trust region on the duration of specific free-time phases is implemented to ensure convergence.
Numerical results are presented to show the effectiveness of the proposed method.
In particular, an analysis of the sensitivity of the achievable payload with respect to different splash-down locations is carried out.
Also, the robustness of the proposed approach with respect to the initialization is studied, along with its computational efficiency, investigating the potential applicability to real-time guidance of future launch vehicles.

\section{Original Problem Formulation}
In this section, the ascent trajectory is first divided into multiple phases to account for different guidance programs, coasting phases, and mass discontinuities. 
Second, the equations of the 3-DoF motion 
of the launch vehicle are derived. 
Finally, the constraints and objective function of the addressed OCP are outlined.


\subsection{Flight Strategy \& Phase Sequence}
A launch vehicle is a system that, from liftoff to payload release, flies through variable conditions and thus requires different guidance programs to meet all mission requirements.
Moreover, the ascent of a multistage rocket consists of a sequence of propelled and coasting arcs, and features the separation of inert masses at each stage burnout.
To effectively tackle these specificities in the optimization process, the corresponding OCP must be cast as a multi-phase problem.

During the first few seconds after liftoff, the rocket has to retain a vertical attitude in order to fly above the launch tower height and safely clear the site.
Then, a programmed rotation maneuver, referred to as \emph{pitch-over}, starts steering the vehicle axis off from its vertical attitude and eventually aligns it with the relative-to-atmosphere velocity.
In the remainder of the atmospheric flight, the rocket is prescribed to keep heading in the direction of the relative wind to minimize the transverse aerodynamic load.
This is called a \emph{zero-lift gravity turn} (ZLGT) maneuver, since it exploits gravity to steer the vehicle while retaining a null angle of attack.
Finally, once the rocket reaches the sufficiently rarefied layers of the atmosphere, an optimal guidance program can be followed.
This usually corresponds to a Hohmann-like maneuver, meaning that the upper stage performs two burns separated by a long coasting arc.


\begin{figure}[h]
    \begin{center}
        \includegraphics[width=\linewidth]{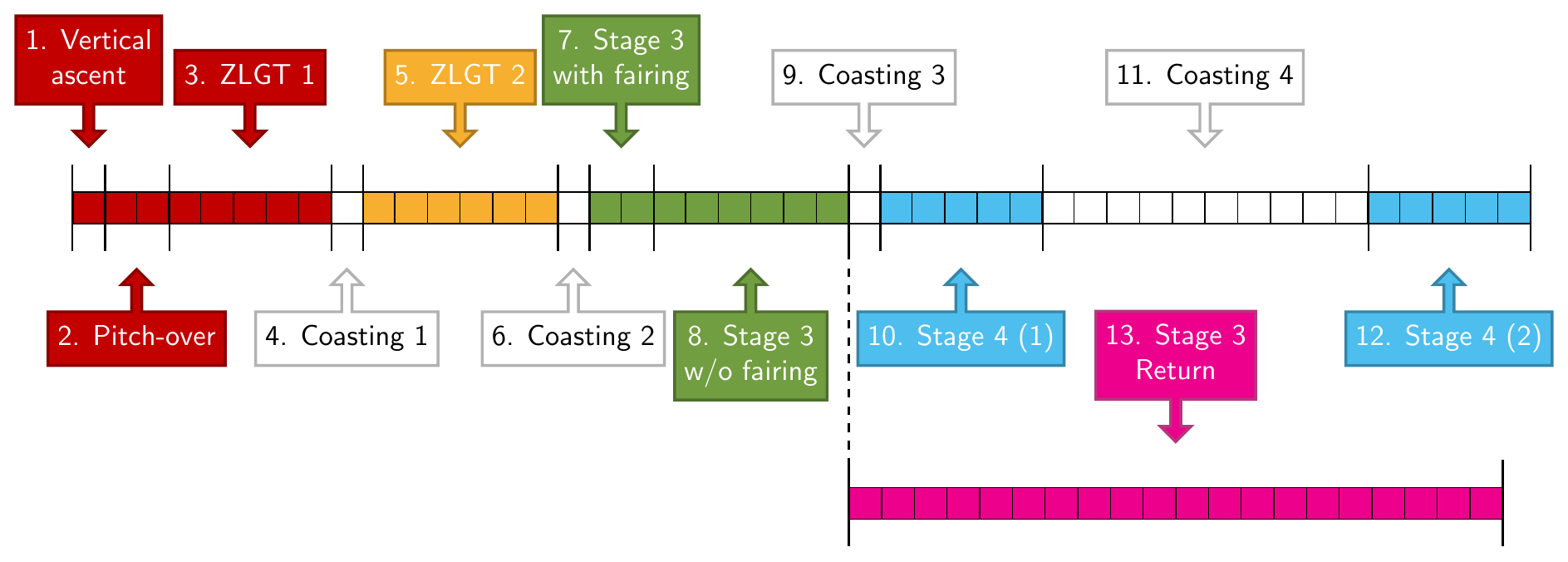}
        \caption{Phases of the optimal control problem}
        \label{fig:phases}
    \end{center}
\end{figure}

The considered phase sequence for VEGA is illustrated schematically in Fig.~\ref{fig:phases} and represents the typical flight strategy of a four-stage launch vehicle.
Note that the phases are numbered progressively from \VerticalAscent{} to \Return{} in chronological order, with the relevant exception of the return phase, which, despite being the \nth{\Return} arc, chronologically starts at the burnout of the third stage, i.e., at the end of phase 8, and takes place concurrently with phases \CoastingThree{}--\StageFourTwo{}.
Hereinafter, let \smash{$t_0^{(i)}$} and \smash{$t_f^{(i)}$} denote the initial and final time of the $i$-th phase.
For the sake of simplicity, if no phase superscript is specified, then $t_0$ and $t_f$ denote the liftoff time $\smash{t_0^{(\VerticalAscent{})}}$ and the fourth stage burnout $\smash{t_f^{(\StageFourTwo{})}}$, respectively.
Likewise, let $t_R$ denote the return time of the spent stage \smash{$t_f^{(\Return)}$}.

The first stage ascent is divided into three phases to properly account for the different guidance programs: vertical ascent (\VerticalAscent{}), pitch-over (\PitchOver{}), and gravity turn (\ZlgtOne{}).
The gravity turn maneuver continues for the entire second stage burn, so phase \ZlgtTwo{} lasts for its whole operation.
The third stage operates at sufficiently high altitudes and can adopt an optimal guidance program as aerodynamic loads do not represent a concern anymore.
Since also the thermal environment is less critical, during the third stage flight, the payload fairing is jettisoned. 
In order to efficiently handle the related mass discontinuity, the third stage operation is split into phases \StageThreeOne{} and \StageThreeTwo{} in correspondence of the jettisoning.
VEGA's last stage performs a Hohmann-like maneuver and its flight is conveniently split into two burn phases, \StageFourOne{} and \StageFourTwo{}, separated by a coasting one, phase \CoastingFour{}.
Note that three other brief coasting arcs (\CoastingOne{}, \CoastingTwo{}, and \CoastingThree{}) are included at each stage separation.
Finally, the return of the third stage is included as phase \Return{} of the OCP.

For the sake of simplicity, we assume a fixed time schedule of phases \VerticalAscent{}--\CoastingThree{}.
Therefore, only the time-lengths of phases \StageFourOne{}--\Return{} are free to be optimized.
Indeed, we assume that the vertical ascent lasts only the few seconds necessary to reach the given clearance altitude over the launchpad and prescribe the pitch-over duration to a value that guarantees that the angle of attack of the launcher is (almost) null at the beginning of the gravity turn.
Moreover, the fairing is supposed to be released after an assigned (small) amount of time 
to guarantee 
the vehicle attitude controllability and the stage full operative conditions at the jettisoning.
The duration of the coasting phases at stage separation are prescribed, with the relevant exception of the Hohmann-like coasting of phase 11.
Finally, the time-lengths of the other time-fixed phases 
are constrained by the (assigned) burn times of each stage.

\subsection{System Dynamics}
The vehicle is modeled as a point mass subject to a 3-DoF translational motion.
Under these assumptions, the state vector $\bm{x}$ is composed of the position vector $\bm{r}$, the velocity vector $\bm{v}$, and the launch vehicle mass $m$:
\begin{equation}
	\bm{x} = \begin{bmatrix}
	x & y & z & v_x & v_y & v_z & m 
	\end{bmatrix}^T
	\label{eq:state_vector}
\end{equation}
Note that the rocket position and velocity are expressed in Cartesian Earth-centered inertial (ECI) coordinates.
In particular, the $x$ axis is in the Earth equatorial plane and passes through the meridian of the launch site at the initial time, the $z$ axis is aligned with Earth angular velocity, and the $y$ axis completes the right-hand frame.
This set of state variables was preferred over the spherical coordinates used in previous works\cite{benedikter2019convexascent,benedikter2019convexrendezvous} because it allows for studying missions toward high inclination orbits without suffering from the numerical issues related to the singularities at the poles. 
As a downside, when using Cartesian coordinates, the terminal conditions result in nonlinear expressions of the state variables.

The launch vehicle is supposed to be subject only to the gravity acceleration $\bm{g}$, the aerodynamic drag $\bm{D}$, and the engine thrust $\bm{T}$.
A Keplerian gravitational model is assumed and the drag force is:
\begin{equation}
    \bm{D} = -\frac{1}{2} C_D S \rho {v}_{\text{rel}} \bm{v}_{\text{rel}}
    \label{eq:drag}
\end{equation}
where $C_D$ is the drag coefficient, assumed to be constant, $S$ is the reference surface, $\rho$ is the atmospheric density, and $\bm{v}_{\text{rel}}$ is the relative-to-atmosphere velocity.
The latter is given by:
\begin{equation}
    \bm{v}_{\text{rel}} = \bm{v} - \bm{\omega}_E \times \bm{r} 
    \label{eq:v_rel}
\end{equation}
where $\bm{\omega}_E$ is Earth's angular velocity vector.

As for the propulsive system, each stage is characterized by a vacuum thrust law $T_{\text{vac}}(t)$ and an ejected mass flow rate $\dot{m}_e(t)$.
However, note that the actual thrust magnitude acting on the system depends also on the external pressure $p$ as:
\begin{equation}
    T = T_{\text{vac}}(t) - p A_e 
\end{equation}
where $A_e$ is the nozzle exit area.

While the thrust magnitude is prescribed by the engine characteristics and the atmospheric conditions, the thrust direction vector $\bm{\hat{T}}$ must be optimized and represents the control $\bm{u}$.
Its elements are expressed in the ECI frame and, since $\bm{\hat{T}}$ is a unit vector, the following relationship must be satisfied:
\begin{equation}
    \hat{T}_x^2 + \hat{T}_y^2 + \hat{T}_z^2 = 1
	\label{eq:thrust_direction_equality_path_con}
\end{equation}

The resulting equations of motion $\bm{\dot{x}} = \bm{f}(\bm{x}, \bm{u}, t)$ are: 
\begin{align}
    \dot{\bm{r}} &= \bm{v} \label{eq:original_ODE_r} \\
    \dot{\bm{v}} &= -\frac{\mu}{r^3} \bm{r} + \frac{T_a}{m} \bm{\hat{T}} + \frac{T_b - D}{m} \bm{\hat{v}}_{\text{rel}} + \frac{T_c}{m} \bm{\hat{r}} \label{eq:original_ODE_v} \\
    \dot{m} &= -\dot{m}_e \label{eq:original_ODE_m}
\end{align}
where the thrust magnitude $T$ is fictitiously split into three contributions to account for the different guidance programs. 
$T_a$ represents the optimally controlled thrust contribution,  while $T_b$ and $T_c$ are always parallel to the relative velocity and to the radial direction, respectively. 
Note that only one of the three terms can be non-null at a given time: 
during the ZLGT arcs only $T_b$ is active; 
similarly, in the vertical ascent phase, $T_c$ is the only active contribution;
finally, $T_a$ is non-zero during the other propelled arcs.

\subsection{Optimal Control Problem}
The goal of the optimization is to determine the control law and other mission parameters, such as the duration of free-time arcs, that maximize the payload mass injected into the target orbit.
In the present work, the propellant and inert masses of the four stages, denoted by $m_{p, i}$ and $m_{\text{dry}, i}$ for $i = 1, \dots, 4$, are supposed to be assigned.
Therefore, we can equivalently decide to maximize the final mass, since it differs from the payload mass by a constant value.
Let the OCP be cast as a minimum problem, then the cost function $J$ to minimize is:
\begin{equation}
    J = -m(\tf)
	\label{eq:objective_m_minimize}
\end{equation}

Besides the payload maximization, the optimization must take into account all mission requirements, which are transcribed as differential, boundary, and path constraints.
The differential constraints are associated with the equations of motion \eqref{eq:original_ODE_r}--\eqref{eq:original_ODE_m}.
The boundary conditions include the initial, terminal, and linkage constraints.
While the initial launcher mass is free to be optimized, the initial position and velocity are completely assigned.
In particular, the launcher initial position corresponds to the launch base location at liftoff $\bm{r}_{\text{LB}}$ and its velocity is equal to the eastward inertial velocity due to Earth rotation:
\begin{align}
    \bm{r}(\ti) &= \bm{r}_{\text{LB}} \label{eq:initial_position} \\
    \bm{v}(\ti) &= \bm{\omega}_E \times \bm{r}_{\text{LB}} \label{eq:initial_velocity}
\end{align} 

The terminal conditions at $\tf$ concern the semi-major axis $a_{\text{des}}$, eccentricity $e_{\text{des}}$, and inclination $i_{\text{des}}$ of the desired orbit.
In particular, assuming a circular target orbit, the following set of constraints can be enforced:
\begin{align}
    x(\tf)^2 + y(\tf)^2 + z(\tf)^2 &= a_{\text{des}}^2 \label{eq:final_radius_nonlinear} \\
    v_{x}(\tf)^2 + v_{y}(\tf)^2 + v_{z}(\tf)^2 &= \mu / a_{\text{des}} \label{eq:final_velocity_nonlinear} \\
    \bm{r}(\tf) \cdot \bm{v}(\tf) &= 0 \label{eq:final_radial_velocity} \\
    x(\tf) v_y(\tf) - y(\tf) v_x(\tf) &= h_{z, \text{des}}
    \label{eq:final_ainc_circular}
\end{align}
Eqs.~\eqref{eq:final_radius_nonlinear} and \eqref{eq:final_velocity_nonlinear} constrain the semi-major axis of the final orbit to be $a_{\text{des}}$.
Eq.~\eqref{eq:final_radial_velocity} guarantees that the radial velocity is zero at payload release, thus, combined with the previous conditions on position and velocity magnitude, ensuring that the final orbit is circular.
Finally, Eq.~\eqref{eq:final_ainc_circular} derives from the expression of the inclination in ECI coordinates:
\begin{equation}
    i = \cos^{-1}\left(\frac{x v_y - y v_x}{h}\right)
    \label{eq:inclination_long}
\end{equation}
Indeed, since the angular momentum $h$ of the target orbit is known and equal to $\sqrt{\mu a_{\text{des}}}$, Eq.~\eqref{eq:inclination_long} can be conveniently expressed as in Eq.~\eqref{eq:final_ainc_circular}, with $h_{z, \text{des}} = \cos i_{\text{des}} \sqrt{\mu a_{\text{des}}}$.

Terminal conditions are prescribed also for the return of the burned-out VEGA's third stage:
\begin{gather}
    x(t_R)^2 + y(t_R)^2 + z(t_R)^2 = R_E^2 \label{eq:final_r_reentry_nonlinear} \\
    z(t_R) = R_E \sin\varphi_{R, \text{des}} \label{eq:final_LAT_reentry}    
\end{gather}
where 
$R_E$ denotes the Earth radius.
Equation~\eqref{eq:final_r_reentry_nonlinear} constrains the final altitude of the returned stage to be null and Eq.~\eqref{eq:final_LAT_reentry} constrains the splash-down location to a given latitude $\varphi_{R, \text{des}}$.
Note that, for missions toward polar or quasi-polar orbits (e.g., Sun-synchronous orbits), constraining the latitude is equivalent to constraining the splash-down distance from the launch base,
since the orbital plane of the trajectory is selected during the pitch-over maneuver and remains (almost) constant in the remainder of the ascent.
This turns out to be a simple, yet effective, way to impose the splash-down constraint, as it consists in assigning just the final value of the $z$ variable,
and it well suits the VEGA target orbits, which are typically high inclination orbits.
Extension to the case of constraining (also) the longitude is straightforward.

Since the problem consists of multiple phases, proper linkage conditions must be enforced at each internal boundary.
All state variables are continuous at boundaries, with the relevant exception of mass, which features a discontinuity at each stage separation:
\begin{align}
    m(t_0^{(\CoastingOne{})}) &= m(t_f^{(\ZlgtOne{})}) - m_{\text{dry}, 1} \label{eq:mass_lkg_con_1} \\
    m(t_0^{(\CoastingTwo{})}) &= m(t_f^{(\ZlgtTwo{})}) - m_{\text{dry}, 2} \label{eq:mass_lkg_con_2} \\
    m(t_0^{(\CoastingThree{})}) &= m(t_f^{(\StageThreeTwo{})}) - m_{\text{dry}, 3} \label{eq:mass_lkg_con_3}
\end{align}
Likewise, at the fairing jettisoning, a mass discontinuity must be accounted for:
\begin{equation}
    m(t_0^{(\StageThreeTwo{})}) = m(t_f^{(\StageThreeOne{})}) - m_{\text{fairing}} \label{eq:mass_hs_lkg_con}
\end{equation}
Finally, the return phase initial boundary corresponds to the third stage burnout, so the following linkage conditions must be enforced:
\begin{align}
    \bm{r}(t_0^{(\Return{})}) &= \bm{r}(t_f^{(\StageThreeTwo{})}) \label{eq:return_lkg_position} \\ 
    \bm{v}(t_0^{(\Return{})}) &= \bm{v}(t_f^{(\StageThreeTwo{})}) \label{eq:return_lkg_velocity} \\ 
    m(t_0^{(\Return{})}) &= m_{\text{dry}, 3} \label{eq:return_mass}
\end{align} 
Note that Eq.~\eqref{eq:return_mass} is not properly a linkage condition, 
since the third stage dry mass is known \emph{a priori}.

As mentioned above, the final stage burn is partitioned between two phases 
(\StageFourOne{} and \StageFourTwo{}).
Since we assumed that all the propellant must be consumed, the sum of the time-lengths of the two firings must be equal to the overall stage burn time $t_{b, 4}$:
\begin{equation}
    \Delta t^{(\StageFourOne{})} + \Delta t^{(\StageFourTwo{})} = t_{b, 4}
    \label{eq:Dt_stage4_constraint}
\end{equation}
where $\Delta t^{(i)} = t_f^{(i)} - t_0^{(i)}$.

Because of the high relative velocity during the atmospheric flight, the rocket undergoes severe thermal conditions.
So, the payload must be protected by the fairing during the initial phases of the ascent.
Nevertheless, once the atmospheric density has decreased enough, the fairing is jettisoned in order to reduce the inert mass as soon as possible.
As a consequence, the payload is directly exposed to the heat flux, which must not exceed a given value.
So, the following path constraint is included in the formulation for phases \StageThreeTwo{}--\StageFourTwo{}:
\begin{equation}
    \dot{Q} = \frac{1}{2} \rho v_{\text{rel}}^3 \leq \dot{Q}_{\text{max}}
    \label{eq:heat_flux_nonlinear}
\end{equation}


\section{Convex Transcription}

In this section, the OCP is formulated as a second-order cone programming (SOCP) problem.
SOCP is a special class of convex programming that is characterized by a linear objective, linear equality constraints, and second-order cone constraints. 
SOCP allows for representing quite complex constraints and can be solved with a small computational effort by means of highly-efficient interior point methods~\cite{alizadeh2003second}.
Since the original problem 
is not convex, it is converted into a SOCP problem via several convexification methods.
First, a convenient change of variables, which produces control-affine dynamics, is proposed.
Second, a control constraint is relaxed into a second-order cone constraint.
The remaining nonconvexities are then tackled via successive linearization.
Virtual controls and buffer zones are introduced to prevent possible artificial infeasibility due to the linearization.
Finally, the continuous-time problem is discretized via a $hp$ pseudospectral method.

\subsection{Change of Variables}
The equations of motion \eqref{eq:original_ODE_r}--\eqref{eq:original_ODE_m} are highly nonlinear in both state and control variables, and thus represent a source of nonconvexity.
A successive linearization of these equations would produce linear constraints, but, due to the coupling of states and controls, high-frequency jitters would show up in the solution process, hindering its convergence\cite{liu2015entry}.
To prevent this undesired behavior, a change of variables is proposed to obtain a control-affine dynamical system.
The new control is introduced:
\begin{equation}
    \bm{u} = \frac{T_a}{m} \bm{\hat{T}} \label{eq:control_vector_cvx}
\end{equation}
Note that $\bm{u}$ includes both the thrust-to-mass ratio $T_{a} / m$ and the thrust direction $\bm{\hat{T}}$.
Replacing the new control in 
Eqs.~\eqref{eq:original_ODE_r}--\eqref{eq:original_ODE_m}
directly produces control-affine equations:
\begin{align}
    \dot{\bm{r}} &= \bm{v} \label{eq:affine_ODE_r} \\
    \dot{\bm{v}} &= -\frac{\mu}{r^3} \bm{r} + \bm{u} 
    + \frac{T_b - D}{m} \bm{\hat{v}}_{\text{rel}} 
    + \frac{T_c}{m} \bm{\hat{r}}
    \label{eq:affine_ODE_v} \\
    \dot{m} &= -\dot{m}_e \label{eq:affine_ODE_m}
\end{align}
So, state and control variables are decoupled and the dynamics can be expressed as:
\begin{equation}
    \bm{f} = \tilde{\bm{f}}(\bm{x}, t) + \tilde{B} \bm{u}
    \label{eq:ODE_split}
\end{equation}
where:
\begin{equation}
    \tilde{B} = \begin{bmatrix}
        \bm{0}_{3 \times 3} \\
        \bm{I}_{3 \times 3} \\
        \bm{0}_{1 \times 3}
    \end{bmatrix}
\end{equation}
$\bm{0}_{m \times n}$ and $\bm{I}_{m \times n}$ denote the null and identity matrix of size $m \times n$.

The new control variables must satisfy Eq.~\eqref{eq:thrust_direction_equality_path_con}, which is reformulated as:
\begin{equation}
    u_x^2 + u_y^2 + u_z^2 = u_N^2
	\label{eq:thrust_direction_equality_path_con_new}
\end{equation}
where the additional variable $u_N$ was introduced:
\begin{equation}
    u_N = \frac{T_a}{m} \label{eq:u_N}
\end{equation}

\subsection{Constraint Relaxation}
\label{subsec:constraint_relaxation}
The path constraint \eqref{eq:thrust_direction_equality_path_con_new} is a nonlinear equality constraint that requires to be convexified in order to be included in the SOCP formulation.
Let us consider its relaxation attained by substituting the equality sign with the inequality sign:
\begin{equation}
    u_x^2 + u_y^2 + u_z^2 \leq u_N^2
	\label{eq:thrust_direction_cone_con}
\end{equation}
Eq.~\eqref{eq:thrust_direction_cone_con} is a convex constraint, in particular a second-order cone constraint.
The inequality sign allows the control variables to be located inside a sphere of radius $u_N$, rather than being constrained on its surface.
Therefore, the convex relaxation defines a larger feasible set than the original one.
Nevertheless, the following proposition ensures that, under mild assumptions, the resulting OCP shares the same solution as the original problem.
Note that the return phase can be temporarily removed from the optimal control problem, as, being an uncontrolled phase, it is not affected by the control constraint relaxation.


\begin{assumption}
    Constraint \eqref{eq:heat_flux_nonlinear} is assumed to be inactive a.e.\footnote{%
    A condition satisfied almost everywhere (a.e.) means that it can be violated only at a finite number of points (a set of measure zero).
}
    in $[t_0, t_f]$.
    \label{assumption:hf}
\end{assumption}

\begin{remark}
    Assumption~\ref{assumption:hf} states that the heat flux constraint is not active over finite intervals of the solution.
    This assumption holds almost always for the ascent problem, since typically the heat flux constraint is active only at isolated points in time, e.g., at the fairing jettisoning.
\end{remark}

\begin{proposition}
Let $\mathcal{P}_{A}$ be the launch vehicle ascent OCP:
\begin{align}
    \mathcal{P}_{A} : \;\;  \min_{\bm{x}, \; \bm{u}, \; t_f} \quad &\eqref{eq:objective_m_minimize} \label{eq:ascent_problem} \\
    \text{s.t.} \quad 
    &\text{
        \eqref{eq:initial_position}--\eqref{eq:final_ainc_circular}, 
        \eqref{eq:mass_lkg_con_1}--\eqref{eq:mass_hs_lkg_con},
        \eqref{eq:Dt_stage4_constraint},
        \eqref{eq:heat_flux_nonlinear},
        \eqref{eq:affine_ODE_r}--\eqref{eq:affine_ODE_m}, 
        \eqref{eq:thrust_direction_equality_path_con_new},
        \eqref{eq:u_N}
        } \nonumber
\end{align}    
Let $\mathcal{P}_{R}$ be the relaxed version of $\mathcal{P}_{A}$ obtained by substituting Eq.~\eqref{eq:thrust_direction_equality_path_con_new} with Eq.~\eqref{eq:thrust_direction_cone_con}, that is:
\begin{align}
    \mathcal{P}_{R} : \;\;  \min_{\bm{x}, \; \bm{u}, \; t_f} \quad &\eqref{eq:objective_m_minimize} \label{eq:relaxed_problem} \\
    \text{s.t.} \quad 
    &\text{
        \eqref{eq:initial_position}--\eqref{eq:final_ainc_circular}, 
        \eqref{eq:mass_lkg_con_1}--\eqref{eq:mass_hs_lkg_con},
        \eqref{eq:Dt_stage4_constraint},
        \eqref{eq:heat_flux_nonlinear},
        \eqref{eq:affine_ODE_r}--\eqref{eq:affine_ODE_m}, 
        \eqref{eq:u_N},
        \eqref{eq:thrust_direction_cone_con}
        } \nonumber
\end{align}
The solution of the relaxed problem $\mathcal{P}_R$ is the same as the solution of $\mathcal{P}_A$. 
That is, if $\{ \bm{x}^\star ; \bm{u}^\star ; t_f^\star \}$ is a solution of $\mathcal{P}_R$, then it is also a solution of $\mathcal{P}_A$ and $u_x^\star(t)^2 + u_y^\star(t)^2 + u_z^\star(t)^2 = u_N^\star(t)^2$ a.e. in $[ t_0 , t_f^\star ]$.
\label{proposition:relaxation}
\end{proposition}
The proof of Proposition~\ref{proposition:relaxation} can be easily obtained by following the same reasoning as in the work by Liu et al.\cite{liu2016exact}, but it is here omitted for the sake of conciseness.
The intuition that motivates Proposition~\ref{proposition:relaxation} lies on the fact that when Eq.~\eqref{eq:thrust_direction_cone_con} is strictly satisfied the engine does not provide the maximum attainable acceleration to the rocket.
Since the goal of the optimization is to maximize the mass injected into a target orbit, it is apparent that such a behavior is suboptimal, and thus will be automatically discarded by the solution procedure.
Finally, note that this relaxation improves the convergence properties of the successive convexification algorithm compared to a linearization of the constraint \eqref{eq:thrust_direction_equality_path_con_new}, since it introduces no approximation and fully preserves the nonlinearity of the original problem.
The benefits of this approach have also been recently investigated and compared to direct linearization by Yang and Liu\cite{yang2019comparison}.

\subsection{Successive Linearization}
Successive linearization is employed to tackle the remaining nonconvexities, which cannot be tackled via lossless convexification.
In particular, the nonconvex constraints are replaced with the first-order Taylor series expansion around a reference solution that is recursively updated.

\subsubsection{Equations of Motion.}
The equations of motion \eqref{eq:affine_ODE_r}--\eqref{eq:affine_ODE_m} are control-affine but still nonlinear in the state variables, thus they must be linearized.
To account for free-time phases in the optimization procedure, we replace time $t$ with $\tau$, a new independent variable defined, for each phase, over a fixed domain $[0, 1]$,
as commonly done in traditional direct methods\cite{betts2010practicalDt}.
Thanks to the unitary domain, the time dilation $\sigma$ between $t$ and $\tau$ is equal to the actual arc time-length:
\begin{equation}
    \sigma = \frac{d t}{d \tau} = t_f - t_0
    \label{eq:dilation_coefficient}
\end{equation}
$\sigma$ is included as an additional optimization parameter for each phase.
The equations of motion are then expressed in terms of $\tau$ and approximated via the first-order Taylor series expansion around a reference solution $\{\bar{\bm{x}}, \bar{\bm{u}}, \bar{\sigma}\}$:
\begin{equation}
    \bm{x}' \vcentcolon= \frac{d \bm{x}}{d \tau} = \sigma \bm{f}\left( \bm{x}, \bm{u}, \tau \right) \approx A \bm{x} + B \bm{u} + \Sigma \sigma + \bm{c}
    \label{eq:linear_ODEs}
\end{equation}
where the following matrices were introduced:
\begin{align}
    A &= \bar{\sigma} \frac{\partial \bm{f}}{\partial \bm{x}} ( \bar{\bm{x}}, \bar{\bm{u}}, \tau ) \label{eq:A_matrix_definition} \\
    B &= \bar{\sigma} \frac{\partial \bm{f}}{\partial \bm{u}} ( \bar{\bm{x}}, \bar{\bm{u}}, \tau ) \label{eq:B_matrix_definition} \\
    \Sigma &= \bm{f} ( \bar{\bm{x}}, \bar{\bm{u}}, \tau ) \label{eq:P_matrix_definition} \\
    \bm{c} &= - (A \bar{\bm{x}} + B \bar{\bm{u}}) \label{eq:C_vector_definition}
\end{align}
Thanks to the change of variables previously carried out, $\bm{f}$ is linear in the control variables, thus the $A$ and $B$ matrices do not depend on the reference solution control $\bar{\bm{u}}$, and $B = \bar{\sigma} \tilde{B}$.
This provides enhanced robustness to the successive linearization sequence as intermediate controls can change significantly among the first iterations\cite{benedikter2019convexrendezvous}.
However, the linearized dynamic equations are still function of the reference controls, but note that, when $\sigma = \bar{\sigma}$, Eq.~\eqref{eq:linear_ODEs} reduces to:
\begin{equation}
    \bm{x}' = A \bm{x} + B \bm{u} + \tilde{\bm{c}}
    \label{eq:linear_ODEs_t_fixed}
\end{equation}
where:
\begin{equation}
    \tilde{\bm{c}} = \bar{\sigma} \tilde{\bm{f}}(\bar{\bm{x}}, \tau) - A \bar{\bm{x}}
    \label{eq:c_tilde}
\end{equation}
For arcs of known duration, 
Eq.~\eqref{eq:linear_ODEs_t_fixed} automatically replaces Eq.~\eqref{eq:linear_ODEs}, but the other arcs may suffer from instability issues when $\sigma$ diverges excessively from the reference value, and some expedient may be necessary to ensure convergence.
In the present application, the return phase does not exhibit any unstable behavior related to $\sigma$, 
but the other free-time phases need further safeguarding constraints on their duration.
In particular, a trust region constraint is imposed on the duration of phases \CoastingFour{} and \StageFourTwo{}:\footnote{Note that phase \StageFourOne{} does not require a trust region as its duration is implicitly related to the one of phase \StageFourTwo{} via Eq.~\eqref{eq:Dt_stage4_constraint}}
\begin{equation}
    | \sigma^{(i)} - \bar{\sigma}^{(i)} | \leq \delta^{(i)} \qquad i = \CoastingFour{}, \StageFourTwo{}
    \label{eq:Dt_trust_region}
\end{equation}
The trust radii $\delta^{(i)}$ are additional optimization variables that are constrained in the interval $[0, \delta_{\text{max}}^{(i)}]$. 
In the authors' experience,
a suitable choice of the upper bound is usually somewhere between 1 and 10\% of $\bar{\sigma}^{(i)}$.
Moreover, to further incentivize $\sigma \approx \bar{\sigma}$, the trust radii are included in the cost function as (slightly) penalized terms by introducing the penalty terms:
\begin{align}
    J_\delta^{(i)} = \lambda_{\delta}^{(i)} \delta^{(i)}
    \qquad i = \CoastingFour{}, \StageFourTwo{}
    \label{eq:delta_penalty_term}
\end{align}
where $\lambda_{\delta}$ are the penalty weights, which should be as small as possible in order not to shadow the originally intended objective and let the optimization autonomously determine the optimal arc time-lengths.

Finally, since the linearization can cause artificial infeasibility\cite{mao2016successive}, a virtual control $\bm{q}$ is included in the dynamics to prevent this undesired phenomenon:
\begin{equation}
    \bm{x}' = A \bm{x} + B \bm{u} + \Sigma \sigma + \bm{c} + \bm{q}
    \label{eq:linear_ODEs_vc}
\end{equation}
The virtual control vector is an unbounded variable that enables to reach any point in the state space in finite time, thus solving the infeasibility issue.
To ensure that its use is limited to otherwise infeasible instances, an additional penalty term is defined:
\begin{equation}
    J_q = \lambda_q P(\bm{q})
    \label{eq:penalty_vc}
\end{equation}
where $\lambda_q$ is the (high) penalty weight and $P(\bm{q})$ a penalty function that we will define upon discretization.

\subsubsection{Boundary Constraints.}
All terminal conditions at payload release \eqref{eq:final_radius_nonlinear}--\eqref{eq:final_ainc_circular} are nonlinear in the state variables and must be linearized as:
\begin{gather}
    \bar{\bm{r}}(t_f) \cdot \bar{\bm{r}}(t_f) 
    + 2 \bar{\bm{r}}(t_f) \cdot (\bm{r}(t_f) - \bar{\bm{r}}(t_f)) 
    = a_{\text{des}}^2 
    \label{eq:final_radius_linearized} \\
    \bar{\bm{v}}(t_f) \cdot \bar{\bm{v}}(t_f) 
    + 2 \bar{\bm{v}}(t_f) \cdot (\bm{v}(t_f) - \bar{\bm{v}}(t_f)) 
    = {\mu}/{a_{\text{des}}} 
    \label{eq:final_velocity_linearized} \\
    \bar{\bm{r}}(t_f) \cdot \bar{\bm{v}}(t_f) 
    + \bar{\bm{v}}(t_f) \cdot (\bm{r}(t_f) - \bar{\bm{r}}(t_f))
    + \bar{\bm{r}}(t_f) \cdot (\bm{v}(t_f) - \bar{\bm{v}}(t_f)) = 0
    \label{eq:final_radial_velocity_linearized} \\
    \bar{v}_y(t_f) (x(t_f) - \bar{x}(t_f)) - \bar{v}_x(t_f) (y(t_f) - \bar{y}(t_f)) - \bar{y}(t_f) v_x(t_f) + \bar{x}(t_f) v_y(t_f) = h_{z, \text{des}}
    \label{eq:final_ainc_linearized_circular}
\end{gather}
Likewise, also the condition on the return final radius \eqref{eq:final_r_reentry_nonlinear} is linearized as:
\begin{equation}
    \bar{\bm{r}}(t_R) \cdot \bar{\bm{r}}(t_R) 
    + 2 \bar{\bm{r}}(t_R) \cdot (\bm{r}(t_R) - \bar{\bm{r}}(t_R)) 
    = R_E^2 
    \label{eq:final_r_reentry_linearized}
\end{equation}
Since also the linearization of the terminal constraints may generate artificial infeasibility, virtual buffer zones are introduced.
In particular, Eqs.~\eqref{eq:final_radius_linearized}--\eqref{eq:final_ainc_linearized_circular} are grouped into a constraint vector $\bm{\chi} = \bm{0}$ and then relaxed as $\bm{\chi} = \bm{w}$, where $\bm{w}$ are free variables, referred to as virtual buffers.
Like virtual control, the virtual buffers should be used only when necessary, so a penalty term is defined:
\begin{equation}
    J_w = \lambda_w \norm{\bm{w}}_1
\end{equation}
where $\lambda_w$ is the (high) penalty weight.

The augmented cost function that includes the trust radii, the virtual control, and the virtual buffer zone penalties is:
\begin{equation}
    J = -m(t_f) + J_\delta^{(\CoastingFour{})} + J_\delta^{(\StageFourTwo{})} + J_q + J_w
    \label{eq:J_augmented}
\end{equation}

\subsubsection{Path Constraints.}
The auxiliary control variable $u_N$ must be equal to the thrust-to-mass ratio at every time and thus
Eq.~\eqref{eq:u_N}
represents a nonlinear path constraint to be linearized as:
\begin{equation}
    u_N = \frac{T_{vac} - p(\bar{\bm{r}}) A_e}{\bar{m}} \left(1 - \frac{m - \bar{m}}{\bar{m}} \right) - 
    \frac{A_e}{\bar{m}} \frac{d p (\bar{\bm{r}})}{d \bm{r}} \cdot (\bm{r} - \bar{\bm{r}})     \label{eq:u_N_path_con_linearized}
\end{equation}
In the same fashion, the linearized heat flux constraint \eqref{eq:heat_flux_nonlinear} is:
\begin{equation}
    \dot{Q}(\bar{\bm{r}}, \bar{\bm{v}}) + \frac{\partial \dot{Q}}{\partial \bm{r}}(\bar{\bm{r}}, \bar{\bm{v}}) \cdot (\bm{r} - \bar{\bm{r}}) + \frac{\partial \dot{Q}}{\partial \bm{v}}(\bar{\bm{r}}, \bar{\bm{v}}) \cdot (\bm{v} - \bar{\bm{v}}) \leq \dot{Q}_{\text{max}}
    \label{eq:heat_flux_linearized}
\end{equation}
where the partial derivatives of the thermal flux with respect to position and velocity are:
\begin{align}
    \frac{\partial\dot{Q}}{\partial \bm{r}} &= \frac{1}{2} \frac{d \rho}{d \bm{r}} v_{\text{rel}}^3 + \frac{3}{2} \rho v_{\text{rel}} \bm{\omega}_E \times \bm{v}_{\text{rel}} \label{eq:dQ_dr} \\
    \frac{\partial\dot{Q}}{\partial \bm{v}} &= \frac{3}{2} \rho v_{\text{rel}} \bm{v}_{\text{rel}} \label{eq:dQ_dv}
\end{align}

\subsection{Discretization}
As a final step, the continuous-time problem must be transcribed into a finite set of variables and constraints to enable the use of numerical algorithms.
In this respect, we employ a $hp$ pseudospectral method.
The $hp$ discretization combines the advantages of $h$ and $p$ schemes, since it exploits the exponential convergence rate of pseudospectral methods in regions where the solution is smooth and introduces mesh nodes near potential discontinuities\cite{darby2011hp}.
Furthermore, compared to $p$ methods, the $hp$ transcription generates sparser problem instances, i.e., with quasi-diagonal matrices, allowing for the use of more efficient numerical routines.

The discretization splits the time domain into multiple subintervals and imposes the differential constraints in each segment via local orthogonal collocation.
In the present paper, we locally employ the Radau pseudospectral method (RPM)\cite{garg2011advances} since it is one of the most accurate and performing pseudospectral methods\cite{garg2010unified}.
The RPM is also a particularly convenient scheme to embed in a $hp$ discretization, as
it avoids redundant control variables at the segment interfaces and provides the optimal control at each mesh point (except for the final node of the final subinterval).
Indeed, the RPM is based on the Legendre-Gauss-Radau (LGR) abscissas,
which include the initial boundary but not the final one.
Locally, this design does not provide the terminal control in each segment, but globally, the ambiguity drops since the final node of a segment corresponds to the initial boundary of the next one, for which, instead, the control is available.



Since details on the implementation of a $hp$ Radau pseudospectral method can be found in the literature\cite{patterson2014gpopsii},
this paper outlines only the major steps of the discretization scheme. 
First, the $hp$ method splits the independent variable domain $\tau \in [0, 1]$ of each phase into $h$ segments by defining a grid $\mathcal{H}$ of $h + 1$ nodes:
\begin{equation}
    0 = \tau_1 < \dots < \tau_{h + 1} = 1
    \label{eq:h_mesh}
\end{equation}
Then, each segment $[\tau_{s}, \tau_{s + 1}]$, 
is discretized as a grid $\mathcal{N}_s$ of $p_s + 1$ nodes:
\begin{equation}
    -1 = \eta_1 < \dots < \eta_{p_s + 1} = 1 
    \label{eq:subinterval_mesh}
\end{equation}
where $p_s$ is the discretization order of the segment and $\eta$ is a new independent variable defined in the interval $[-1, 1]$, which can be mapped to the original domain by the following transformation:
\begin{equation}
    \tau = \frac{\tau_{s + 1} - \tau_{s}}{2}\eta + \frac{\tau_{s + 1} + \tau_{s}}{2}
    \label{eq:eta2tau}
\end{equation}
Since we employ the RPM, the first $p_s$ nodes of each segment correspond to the set of $p_s$ LGR roots and constitute the collocation points $\mathcal{K}_s$.
Note that $\mathcal{K}_s$ is a subset of $\mathcal{N}_s$ since it does not incorporate the terminal boundary $\eta = 1$.

Once the grid is set up, the state and control are discretized over it, and a finite set of variables $(\bm{x}_j^s, \bm{u}_j^s)$ is obtained.
The superscript $s$ denotes the $s$-th segment, while the subscript $j$ refers to the $j$-th node of the segment.
In particular, in each segment, the state is discretized over the set $\mathcal{N}_s$ and approximated using a basis of Lagrange polynomials.
Note that since the state is continuous among the segments of a phase, in the algorithm implementation the same variable is used for both $\bm{x}_{p_s + 1}^s$ and $\bm{x}_{1}^{s + 1}$.
Instead, the control is discretized only at the collocation points $\mathcal{K}_s$, so Lagrange polynomials of degree $p_s - 1$ are used for the approximation.
The final control of the final segment is not included in the discrete problem and it is simply extrapolated from the polynomial approximation of the control signal.



Path constraints are converted into a finite set of algebraic constraints by imposing them at every node, while boundary conditions are imposed only at the initial or final point of $\mathcal{H}$.
To take into account the system dynamics \eqref{eq:linear_ODEs_vc}, the time derivative of the state interpolating polynomial is constrained to be equal to the equations of motion at the collocation points of each segment $s = 1, \dots, h$:
\begin{equation}
    \sum_{j = 1}^{p_s + 1} D_{ij}^s \bm{x}_j^s = \frac{\tau_{s+1} - \tau_{s}}{2} (A_i^s \bm{x}_i^s + B_i^s \bm{u}_i^s + \Sigma_i^s \sigma + \bm{c}_i^s + \bm{q}_i^s) 
    \qquad 
    i = 1, \dots, p_s 
    \label{eq:collocation_radau}
\end{equation}
where the same notation used for discrete-time variables was used for the linearization matrices \eqref{eq:A_matrix_definition}--\eqref{eq:C_vector_definition}.
In Eq.~\eqref{eq:collocation_radau}, $D^s$ denotes the LGR differentiation matrix\cite{garg2010unified}, which can be efficiently computed via barycentric Lagrange interpolation\cite{berrut2004barycentric}.
Finally, similarly to the other continuous-time variables, also the virtual control is discretized over the mesh as $\bm{q}_j^s$.
The resulting set of variables is grouped into a vector $\tilde{\bm{q}}$ and the penalty term introduced in Eq.~\eqref{eq:penalty_vc} can be transcribed as:
\begin{equation}
    J_q = \lambda_q \norm{\tilde{\bm{q}}}_1 
\end{equation}

As a final remark, in this paper, no automatic mesh refinement is implemented.
So, a sufficiently dense grid must be devised \emph{a priori} according to the desired discretization accuracy.

\section{Reference Solution}

The convexification of the original problem nonlinear dynamics and constraints exploits successive linearization,
which replaces the original expressions with a first-order Taylor series expansion around a reference solution $\{\bar{\bm{x}}, \bar{\bm{u}}, \bar{\sigma}\}$.
This section focuses on the reference solution.
First, we outline a simple procedure to design a starting trajectory that allows for convergence.
Then, an improved method for updating the reference solution based on multiple previous iterations is proposed.

\subsection{Initialization}
Sensitivity to the initialization is a major downside of traditional optimization methods.
For instance, indirect methods can achieve convergence only if an accurate first guess is provided.
This is a cumbersome drawback as an initialization is not required only for the trajectory but also for the costate 
and the structure of the constrained arcs, which can be difficult to supply.
On the other hand, direct methods exhibit greater robustness to the initialization, but the discretization of highly-sensitive nonconvex OCPs, such as the one at hand, produces a NLP problem whose solution depends significantly on the first guess.
These limitations motivate the upstream effort put into the careful convexification process.
Indeed, a greater robustness is observed in the devised algorithm compared to traditional direct optimization methods.
Moreover, also compared to our previous work on convex optimization of the ascent problem\cite{benedikter2019convexascent},
the present algorithm shows an enhanced robustness.
The reason for this improvement is due to the $hp$ pseudospectral discretization,
which accurately describes the dynamics and yet retains a sparse problem structure.

The standard way of dealing with the problem at hand is, first, solving the ascent problem without the splash-down constraint, then simulate the return of the spent stages and, if necessary, constrain the splash-down to a safe location. 
In fact, the concern on the splash-down of the spent stages exists only if the simulation of the return trajectory corresponds to an unsafe impact location.
So, phase \Return{} and the related constraints can be omitted at first and focus can be placed on designing a reference solution for phases \VerticalAscent{}--\StageFourTwo{} only.


The present algorithm does not require an accurate initialization, but, rather, in the authors' experience, 
any starting trajectory with an altitude profile always above sea level is sufficient to achieve convergence.
Such trajectories can be easily generated via numerical integration of the original rocket equations of motion \eqref{eq:original_ODE_r}--\eqref{eq:original_ODE_m}.
To set up the forward propagation, the unknown control history, the duration of free-time arcs, and the initial mass must be prescribed.
In general, designing the control laws may be a complex task, but if the atmosphere is removed from the dynamics and by choosing a small value of $m(t_0)$, i.e., which corresponds to a small payload mass, even trivial control laws can produce acceptable trajectories.

During the pitch-over, the elevation $\phi$, i.e., the angle between the thrust direction and the local horizontal,
is prescribed to vary linearly from \SI{90}{\degree} to a final value,
commonly referred to as \emph{kick} angle,
here denoted as \smash{$\phi(t_f^{(\PitchOver{})})$}. 
Phase~\PitchOver{} is also assumed to take place in an fixed inertial plane; therefore, the thrust azimuth $\psi$, i.e., the angle measured clockwise from the north direction to the thrust vector, is kept equal to a constant value \smash{$\psi^{(\PitchOver)}$}.
While the kick angle must be guessed, a systematic way of choosing \smash{$\psi^{(\PitchOver)}$} is selecting the value that, under the non-rotating Earth assumption, allows for reaching the target orbit plane without further out-of-plane maneuvers:
\begin{equation}
    \psi^{(\PitchOver)} = \sin^{-1}\left(\frac{\cos(i_{\text{des}})}{\cos(\varphi_{LB})}\right)
    \label{eq:pitch_over_azimuth}
\end{equation}
where $\varphi_{LB}$ is the latitude of the launch base.
Unfortunately, when $i_{\text{des}} < \varphi_{LB}$ the previous formula does not hold anymore and an \emph{ad hoc} value must be provided for \smash{$\psi^{(\PitchOver)}$}.
For stages 3 and 4, even simpler control laws can be devised.
Indeed, the tentative solutions are designed such that
the orbital plane is kept constant after the second stage burnout; 
so, the thrust vector is constrained in the orbital plane.
The only control to prescribe is the elevation angle $\phi$, which is kept null for the entire operation of both final stages.
Finally, the subdivision of the fourth stage burn and the duration of the intermediate coasting must be chosen.

To sum up, the only variables necessary for generating a tentative solution are: (i) the initial mass $m(t_0)$ (or, equivalently, the payload mass),
(ii) the kick angle \smash{$\phi(t_f^{(\PitchOver{})})$},
and (iii) the time-lengths of phases \StageFourOne{}--\StageFourTwo{}.
These values should be set on the basis of the specific launch vehicle and target orbit.
Nevertheless, 
their choice does not represent an arduous task, since a wide range of values can generate acceptable trajectories.

\subsection{Recursive Update}

At every iteration of the successive convexification algorithm, the reference solution must be updated.
Traditional successive linearization algorithms solve the $i$-th SOCP problem by linearizing the constraints around the $(i - 1)$-th solution.
Instead, we employ an improved method, named \emph{filtering}\cite{benedikter2019convexrendezvous}, which consists in computing the reference solution for the $i$-th SOCP problem as a weighted sum of the $K$ previous solutions:
\begin{equation}
	\bar{x}^{(i)} = \sum_{k = 1}^{K} \alpha_k x^{\max \{0, (i - k)\} }
\end{equation}
where $\alpha_k$ are constant weights and $x^{(i)}$ denotes the solution to the $i$-th subproblem.
Note that if $i < K$ then the initial reference solution $x^{(0)}$ appears multiple times in the sum.

The proposed technique
adds another layer of algorithmic robustness to the successive convexification procedure.
In fact, it has been observed that sequences solved with $K = 1$ suffered from instability issues, mainly related to artificial unboundedness.
The common approach to unboundedness is adding a trust region constraint that limits the search space to the neighborhood of the reference solution\cite{mao2016successive}.
However, if the reference solution is far from the optimal one, constraining the search space may cause convergence toward a suboptimal solution.
Instead, filtering efficiently solves the unboundedness issue and it does not affect the optimality of the attained solution, as no additional constraint or penalty term is included in the SOCP formulation.

The same parameters used in a different application\cite{benedikter2019convexrendezvous} revealed to be effective also for the problem at hand and thus are deemed as the most performing combination. 
Specifically, the three last solutions are used ($K = 3$) and the corresponding weights are reported in Table~\ref{tab:weights_filtering}.
Nevertheless, also different values of $K$ and of the weights can achieve convergence in a wide range of missions.
Eventually, the sequential algorithm terminates when the difference between the computed solution and the reference one converges below an assigned tolerance:
\begin{equation}
    \norm{\bm{x} - \bar{\bm{x}}}_\infty < \epsilon_{\text{tol}}
    \label{eq:successive_cvx_termination_condition}
\end{equation}

\begin{table}[h]
    \caption{Filtering weights}
    \label{tab:weights_filtering}
    \centering
    \begin{tabular}{c c c}
        \hline
        $\alpha_1$ & $\alpha_2$ & $\alpha_3$ \\
        \hline
        6/11 & 3/11 & 2/11 \\
        \hline
    \end{tabular}
\end{table}

\section{Numerical Results}

In this section, numerical results are presented to show the effectiveness of the proposed approach. 
The described algorithm has been implemented in C++ 
using Gurobi\cite{gurobi} as SOCP solver.
The values of the penalty weights and convergence tolerance are reported in Table~\ref{tab:algoritm_parameters}.
Also, since scaling is key to the effectiveness of any numerical algorithm, we take the Earth radius, the corresponding circular orbit velocity, and a reference mass of \SI{10000}{\kilo\g} as normalization factors.

\begin{table}[h]
    \centering
    \caption{Algorithm parameters}
    \label{tab:algoritm_parameters}
    \begin{tabular}{c c c c c c}
    \hline
    \bf Parameter & $\lambda_\delta^{(\CoastingFour{})}$ & $\lambda_\delta^{(\StageFourTwo{})}$ & $\lambda_q$ & $\lambda_w$ & $\epsilon_{\text{tol}}$ \\
    \hline
    \bf Value & \num{1e-4} & \num{1e-4} & \num{1e4} & \num{1e4} & \num{1e-4} \\
    \hline
    \end{tabular}
\end{table}

The data used to model the VEGA launch vehicle are summarized in Table~\ref{tab:rocket_data}.
The main assumption concerns the thrust and mass flow rate history of the stages, which are approximated as linear functions of time.
Nevertheless, the total impulse is retained and the other quantities are quite accurate, so the overall model is representative of the real system performance.
Other design values include the fairing mass $m_{\text{fairing}}$ (\SI{535.3}{\kilo\g}), 
the drag coefficient $C_D$ (\num{0.381}), 
and the reference surface $S$ (\SI{9.079}{\square\m}).
Although a realistic aerodynamic model would be needed to accurately predict the splash-down location, for this work, in a simplified manner, the same coefficients are used also for the return phase. 
Notwithstanding, the algorithm can consider more realistic aerodynamic characterizations of the launch vehicle and stage return.
Finally, the U.S. Standard Atmosphere 1976 model is used to evaluate the air density and pressure as functions of the altitude\cite{us1976atm}.


\begin{table}[h]
    \centering
    \caption{VEGA-like rocket data}
    \label{tab:rocket_data}
    \begin{tabular}{l c c c c c}
    \hline
    & \textbf{Stage 1} & \textbf{Stage 2} & \textbf{Stage 3} & \textbf{Stage 4} & \textbf{Unit} \\
    \hline
    $m_p$ & \num{87898} & \num{23926} & \num{10006} & \num{397.6} & \si{\kilo\g} \\
    $m_{\text{dry}}$ & \num{8417.7} & \num{2563.8} & \num{1326.5} & \num{813.7} & \si{\kilo\g} \\
    $t_b$ & \num{102.0} & \num{75.0} & \num{110.0} & \num{502.1} & \si{\s} \\
    $T_{\text{vac}}(0)$ & \num{2827.37} & \num{1075.73} & \num{299.81} & \num{2.4509} & \si{\kilo\N} \\
    $T_{\text{vac}}(t_b)$ & \num{1884.91} & \num{717.15} & \num{221.60} & \num{2.4509} & \si{\kilo\N} \\
    $\dot{m}_e(0)$ & \num{1034.09} & \num{382.81} & \num{104.61} & \num{0.7919} & \si{\kilo\g\per\s} \\
    $\dot{m}_e(t_b)$ & \num{689.40} & \num{255.21} & \num{77.32} & \num{0.7919} & \si{\kilo\g\per\s} \\
    $A_e$ & \num{3.092} & \num{1.697} & \num{1.183} & \num{0.07} & \si{\square\m} \\
    \hline
    \end{tabular}
\end{table}


As for the discretization, Table~\ref{tab:mesh} reports $h$ and $p$ for every phase.
The values have been devised in a heuristic way in order to meet the desired discretization accuracy.
In particular, 
since phases \VerticalAscent{}--\StageFourTwo{} are relatively brief and do not feature rapidly changing dynamics, no internal subdivision is necessary and $h$ is simply set to 1.
Instead, since a high number of nodes are required to capture the reentry dynamics and high-order approximating polynomials suffer from numerical issues, the return phase is split into 10 equally-spaced segments.
In each segment, the same discretization order $p$ is used.

\begin{table}[h]
    \centering
    \caption{Discretization segments, order, and nodes in each phase}
    \label{tab:mesh}
    \medskip
    \begin{tabular}{l c c c c c c c c c c c c c}
    \hline
    \bf Phase & 1 & 2 & 3 & 4 & 5 & 6 & 7 & 8 & 9 & 10 & 11 & 12 & 13 \\
    \hline
    $h$ & 1 & 1 & 1 & 1 & 1 & 1 & 1 & 1 & 1 & 1 & 1 & 1 & 10 \\
    $p$ & 5 & 5 & 17 & 5 & 19 & 14 & 5 & 19 & 9 & 19 & 19 & 19 & 10 \\
    Nodes & 6 & 6 & 18 & 6 & 20 & 15 & 6 & 20 & 10 & 20 & 20 & 20 & 101 \\
    \hline
    \end{tabular}
\end{table}

The considered case study is a mission toward a \SI{700}{\kilo\m} circular polar Earth orbit ($i_{\text{des}} = \SI{90}{\degree}$).
The vehicle is assumed to take off from the equator in correspondence of the Guiana Space Center meridian.
The time-lengths of the arcs of fixed duration are reported in Table~\ref{tab:time-lengths}.
The threshold on the bearable heat flux is set to \SI{900}{\W\per\square\m}.
First, the optimal ascent trajectory is found, neglecting the splash-down location of the spent stages.
Then, the return phase is included in the OCP and an analysis of the sensitivity of the system performance to different impact points is presented.


\begin{table}[h]
    \centering
    \caption{Time-lengths of time-fixed arcs}
    \label{tab:time-lengths}
    \medskip
    \begin{tabular}{l c c c c c c c c c}
    \hline
    \bf Phase & 1 & 2 & 3 & 4 & 5 & 6 & 7 & 8 & 9 \\
    \hline
    \bf $\Delta t$ (\si{\s}) & \num{4.1} & \num{6.6} & \num{91.3} & \num{6.6} & \num{75.0} & \num{37.3} & \num{5.4} & \num{104.6} & \num{15.4} \\
    \hline
    \end{tabular}
\end{table}

\subsection{Unconstrained Return}

To set up the optimization, a starting reference solution must be provided.
This is generated as described in the previous section and the used parameters are reported in Table~\ref{tab:first_guess_values}.
Note that a very small payload mass $m_{\text{pl}}$ was picked (less than 10\% of the expected optimum) and that the duration of phase \StageFourOne{} is omitted as it can be automatically derived from Eq.~\eqref{eq:Dt_stage4_constraint}.

\begin{table}[h]
    \centering
    \caption{Values used for the generation of the first guess trajectory}
    \label{tab:first_guess_values}
    \begin{tabular}{l c c}
    \hline
    \bf Quantity & \bf Value & \bf Unit \\
    \hline
    $m_{\text{pl}}$ & \num{100.0} & \si{\kilo\g} \\
    $\phi_{\text{ko}}$ & \num{80.0} & deg \\
    $\Delta t^{(\CoastingFour{})}$ & \num{2500.0} & \si{\s} \\
    $\Delta t^{(\StageFourTwo{})}$ & \num{200.0} & \si{\s} \\
    \hline
    \end{tabular}
\end{table}

\begin{figure}[h]
	\centering\includegraphics[width=0.7\linewidth]{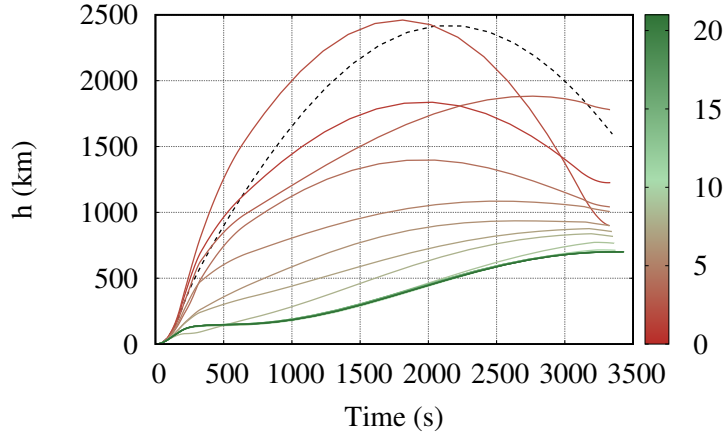}
	\caption{Iteration sequence starting from the initial reference solution (dashed black line). The intermediate trajectories transition from red ($\text{iter} = 1$) to green ($\text{iter} = 22$).}
	\label{fig:iterations}
\end{figure}

The convergence behavior is illustrated in Fig.~\ref{fig:iterations}.
Starting from the initial guess (dashed black line), the intermediate solutions, whose color transitions from red to green, gradually converge to the final trajectory.
Thus, despite the initial reference solution is far from the solution of the OCP, the termination condition \eqref{eq:successive_cvx_termination_condition} is eventually met in 22 iterations.
Note that the virtual buffer zones introduced to relax the terminal constraints are actively exploited in the first 10 iterations.
Indeed, the intermediate subproblems would otherwise be infeasible (even with virtual controls), thus virtual buffers are essential to ensure the recursive feasibility of the sequential process.
Without any specific code optimization, the overall computational time is \SI{12.8}{\s}, so each iteration requires \SI{0.58}{\s} on average.\footnote{%
    The algorithm was tested on a computer equipped with Intel\textregistered{} Core\texttrademark{} i7-7700HQ CPU @ \SI{2.80}{\giga\Hz}.
}
By using a custom SOCP solver and running on dedicated hardware, a further speed-up can be expected, thus enabling potential suitability for real-time guidance.
Moreover, in real-time applications a much more accurate initialization is used, as a nominal trajectory is already available, so fewer iterations are needed, greatly reducing the computational burden.

\begin{figure}[h]
	\centering\includegraphics[width=0.7\linewidth,trim={4cm 0.5cm 2cm 4.5cm},clip]{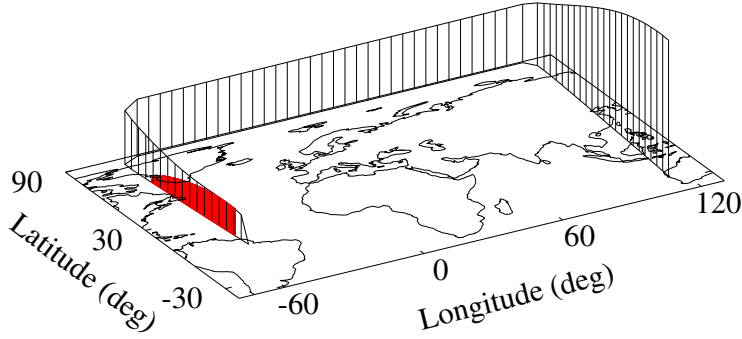}
	\caption{Visualization of the unconstrained trajectory. The third stage return flight is colored in red.}
	\label{fig:groundtrack_3d}
\end{figure}

The unconstrained trajectory is illustrated in 
Fig.~\ref{fig:groundtrack_3d}. 
The figure also reports the simulation of the return phase, which provides the optimal splash-down latitude
($\varphi_{R}^\star = \SI{65.79}{\degree}$).
The accuracy of the converged solution is verified by forward propagation of the original equations of motion \eqref{eq:original_ODE_r}--\eqref{eq:original_ODE_m} using the optimal control laws.
In particular, the discrepancies in the terminal conditions are inspected.
The largest inaccuracy concerns the semi-major axis, but the relative error is below 1\%, corresponding approximately to \SI{100}{\m}, which is in agreement with the finite precision of the SOCP solver.

To validate the quality of attained results, the same problem was solved also using EOS\cite{federici2020eos}, a direct shooting algorithm based on differential evolution that was already successfully employed to solve a similar instance of the problem at hand\cite{Federici2019Integrated}.
The comparison between the two solutions is reported in Table~\ref{tab:comparison_EOS}.
The payload mass difference is approximately \SI{5}{\kilo\g} and is due to the different sets of time-lengths found. 
Indeed, the problem features many local optima with different times but similar costs, so the optimization can converge unpredictably toward one of these.
Nevertheless, note that the difference in cost is minimal, so both solution are acceptable for any practical purpose.
Compared to the convex approach, the main drawback of EOS is the large computational effort required (approximately 20 minutes on the same hardware).

\begin{table}[h]

\sisetup{round-mode = places,round-precision = 2}
\sisetup{round-minimum = 0.01}
\caption{Unconstrained solution compared with the EOS solution.}
\label{tab:comparison_EOS}
\centering 

    \begin{tabular}{l c c c c}
    \hline
    & \bf Convex & \bf EOS & \bf Unit & \bf Variation (\%) \\
    \hline
    $m_{\text{pl}}$ & \num{1400.73} & \num{1396.74} & \si{\kilo\g} & \num{0.2855} \\ 
    $\Delta t^{(\StageFourOne{})}$ & \num{359.71} & \num{357.60} & \si{\s} & \num{0.5892} \\ 
    $\Delta t^{(\CoastingFour{})}$ & \num{2583.50} & \num{2660.70} & \si{\s} & \num{2.9016} \\ 
    $\Delta t^{(\StageFourTwo{})}$ & \num{142.39} & \num{144.50} & \si{\s} & \num{1.4581} \\ 
    \hline
    \end{tabular}
\end{table}
    

\subsection{Parametric Analysis of the Splash-Down Constraint}

Once the optimal solution is obtained, we can investigate the effect of the splash-down constraint.
Therefore, the return phase is added to the OCP along with the corresponding constraints and the impact point is gradually moved from its optimal location to different latitudes.
So, a series of problems with different $\varphi_{R, \text{des}}$ is solved.
Each OCP uses the converged solution of the previous one as initialization.
Since this initialization is quite accurate, on average only \num{8} iterations are required to meet the convergence criterion. 
However, the inclusion of the return phase significantly increases the problem dimension, so each iteration is computationally more demanding (\SI{0.95}{\s} on average).
Nevertheless, the overall process requires a mean computational time of \SI{7.75}{\s}.

\begin{figure}[htb]
    \centering

    \begin{minipage}{0.45\linewidth}
        \centering
        \includegraphics[width=\linewidth]{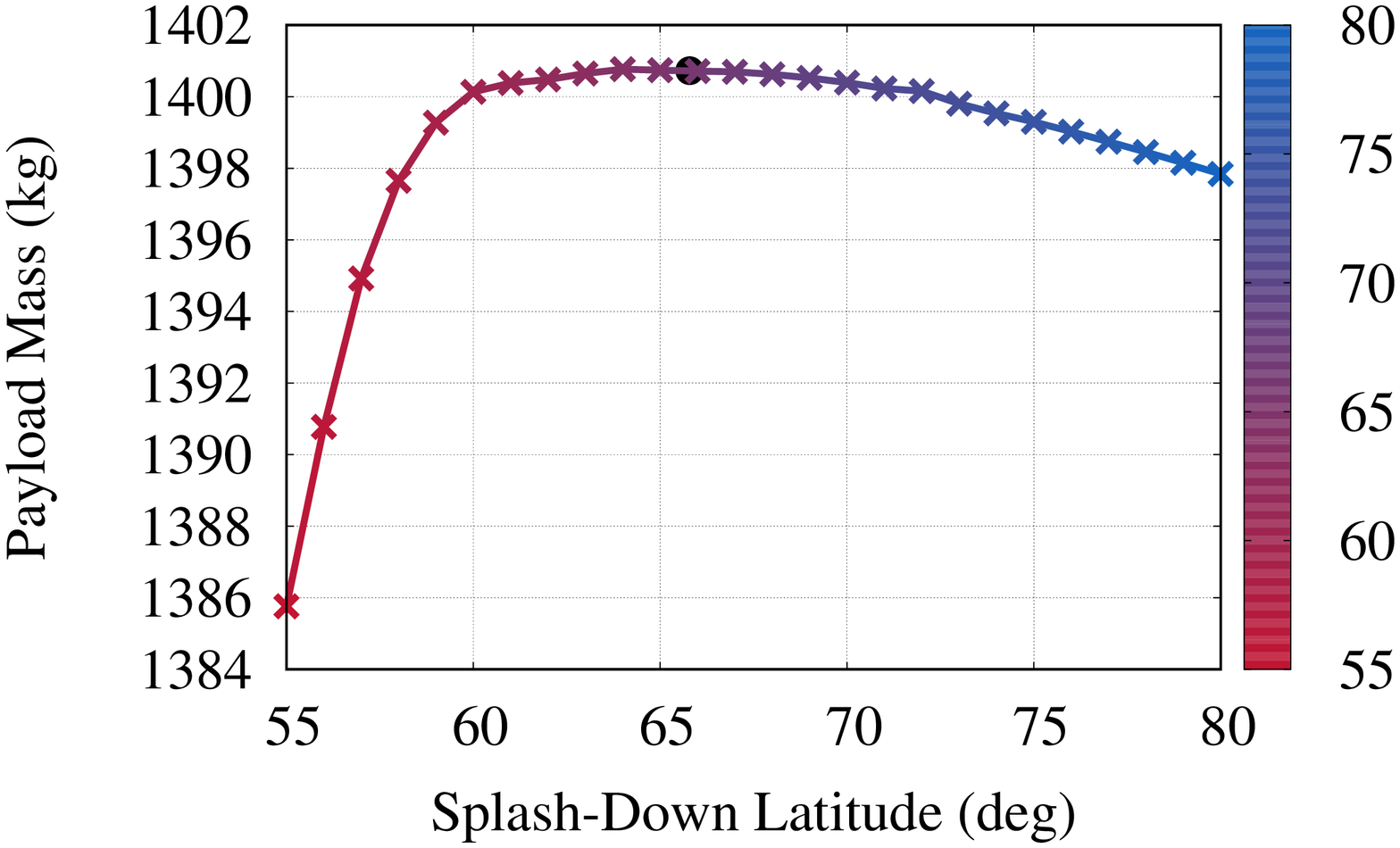}
        \subcaption{Payload mass vs. splash-down latitude}
        \label{fig:payload_lat}
    \end{minipage}
    \hspace{6mm}
    \begin{minipage}{0.45\linewidth}
        \centering
        \includegraphics[width=\linewidth]{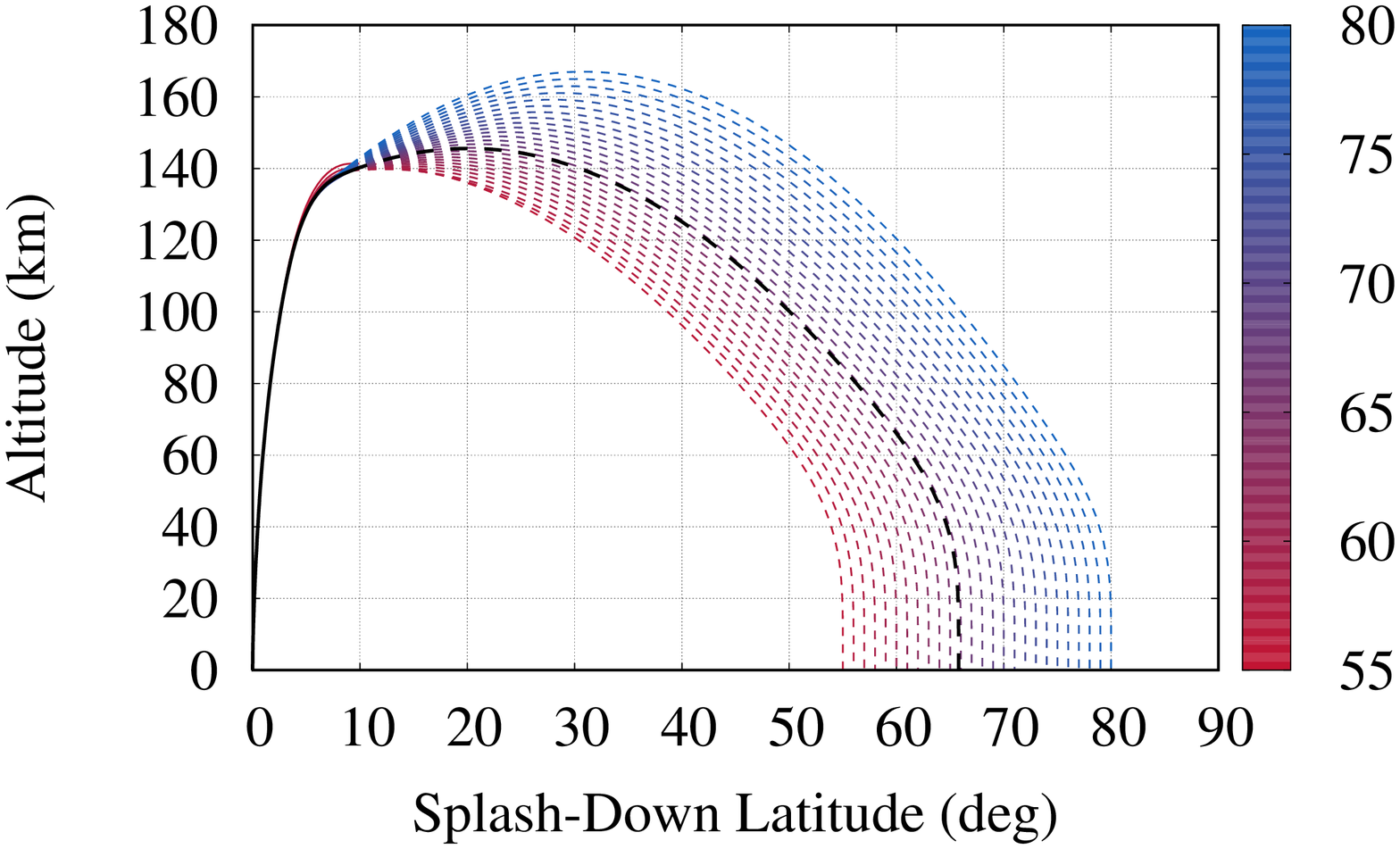}
        \subcaption{Stage 3 trajectories from liftoff to splash-down}
        \label{fig:traj_stage3}
    \end{minipage}

    \begin{minipage}{0.45\linewidth}
        \centering
        \includegraphics[width=\linewidth]{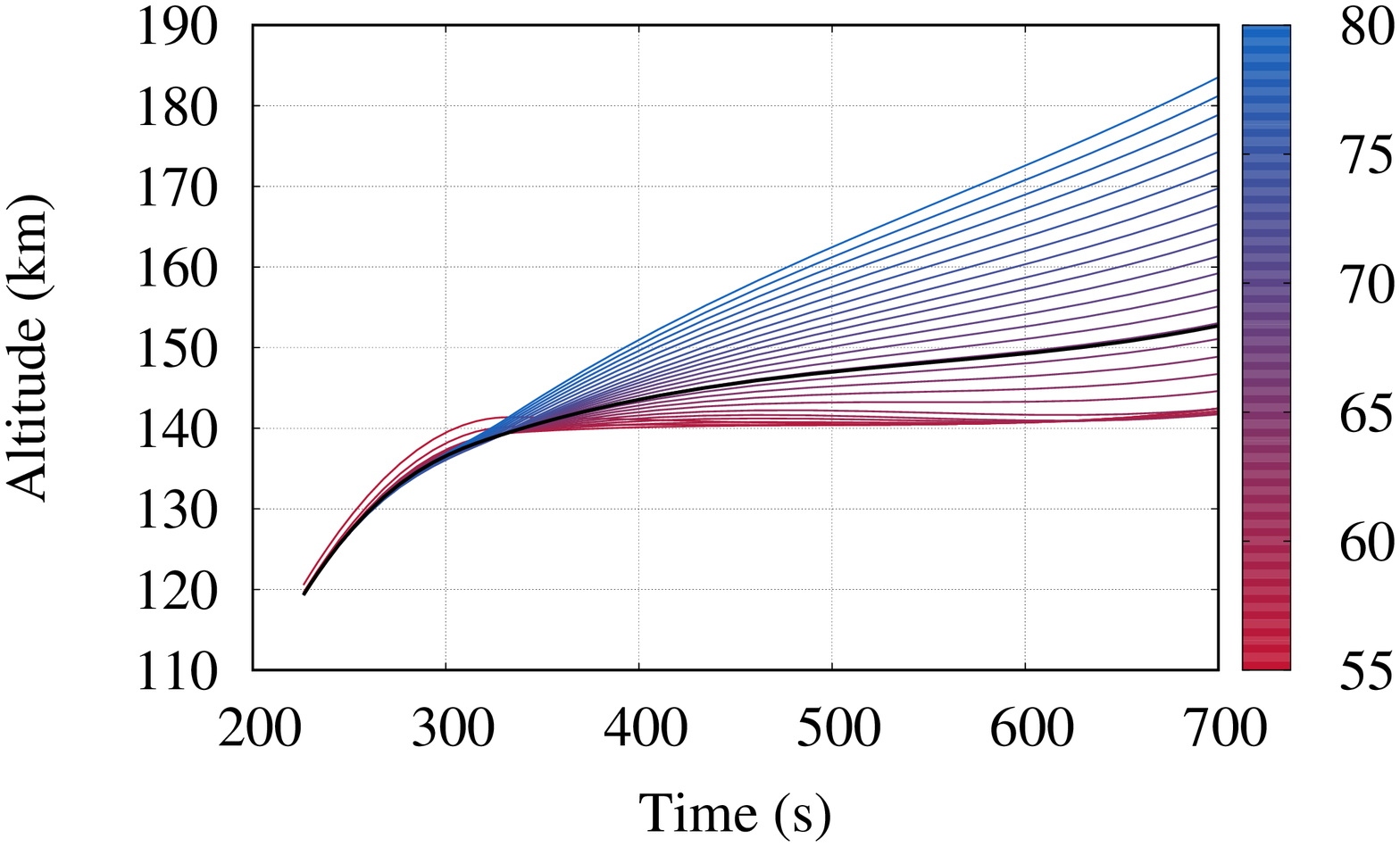}
        \subcaption{Altitude profiles from fairing jettisoning to $t_f^{(\StageFourOne{})}$}
        \label{fig:altitude_HS}
    \end{minipage}
    \hspace{6mm}
    \begin{minipage}{0.45\linewidth}
        \centering
        \includegraphics[width=\linewidth]{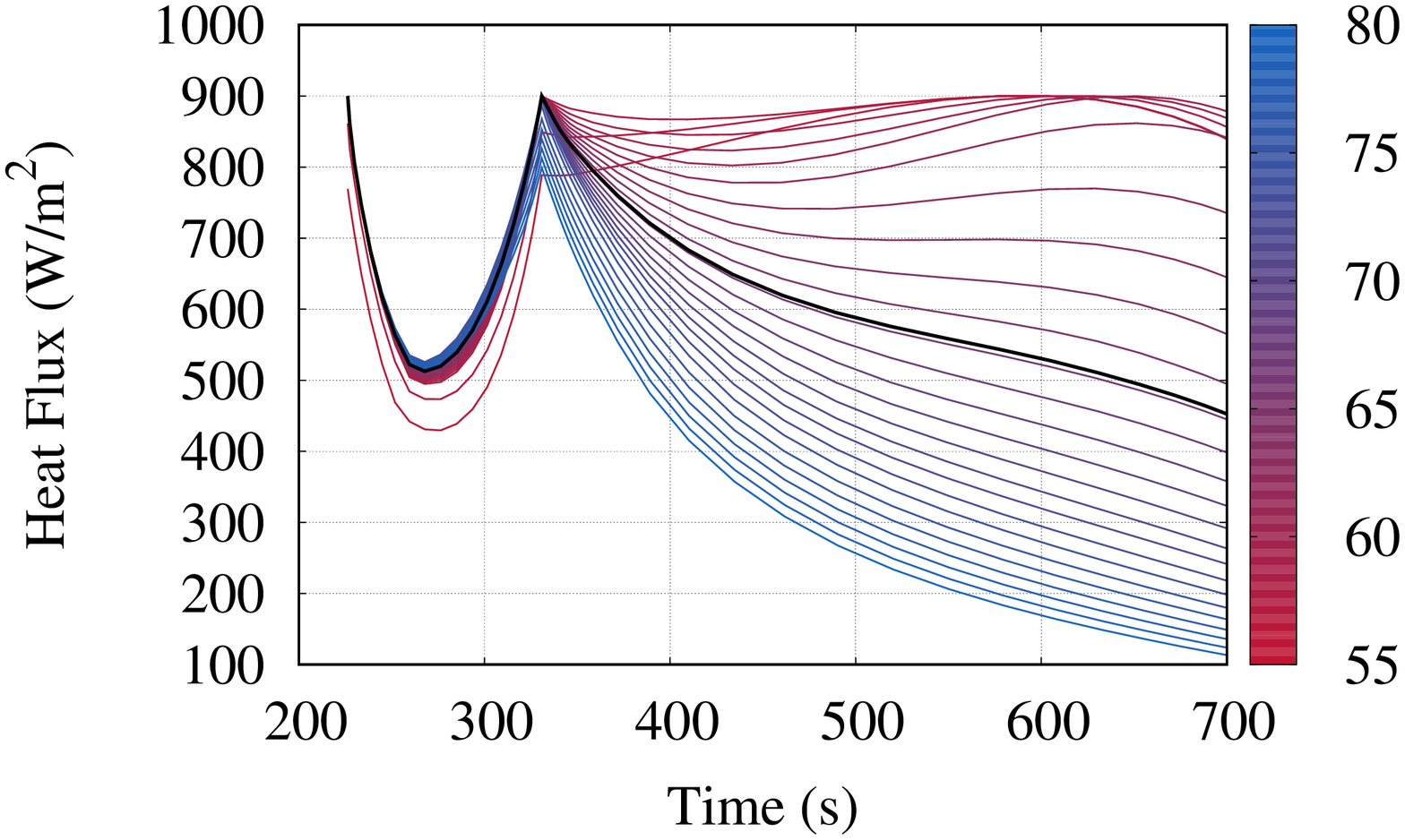}
        \subcaption{Heat flux profiles from fairing jettisoning to $t_f^{(\StageFourOne{})}$}
        \label{fig:hf}
    \end{minipage}
    \caption{Performance, trajectories, and thermal conditions corresponding to different splash-down latitudes. The black curves correspond to the unconstrained solution.}
    \label{fig:splash_down_profiles}
\end{figure}

The optimal payload mass is plotted in Fig.~\ref{fig:payload_lat} as a function of the splash-down latitude.
While moving the spent stage return point significantly changes the trajectory, as shown in Figs.~\ref{fig:traj_stage3} and \ref{fig:altitude_HS}, it does not necessarily hinder the performance.
Indeed, the payload curve is essentially flat in the interval $\varphi_{R} \in [\SI{60}{\degree}, \SI{72}{\degree}]$ and variations only below \SI{1}{\kilo\g} are observed.
When the splash-down location is moved beyond \SI{72}{\degree}, the decrease in performance is more evident, but it still does not represent a concern as a shift of \SI{15}{\degree} causes a loss of only \SI{3}{\kilo\g}.
Instead, moving the splash-down point closer than \SI{60}{\degree} appears more critical, as a greater performance drop is observed.
Nevertheless, even constraining the third stage to fall \SI{10}{\degree} closer than $\varphi_{R}^\star$ results in a payload reduction by only 1\% of its optimal value. 

It is worth studying how the heat flux and splash-down constraints interact with each other.
Fig.~\ref{fig:hf} shows the heat flux history that the payload undergoes from the fairing jettisoning until the end of the first firing of stage 4.
Red curves correspond to splash-down locations at lower latitudes, i.e., closer to the launch site, while blue ones are associated with high-latitude returns.
In all trajectories for which $\varphi_{R} \geq \SI{57}{\degree}$, the heat flux constraint is active only at the boundaries of phase \StageThreeTwo{}, so Assumption~\ref{assumption:hf} holds in all these cases.
Instead, the heat flux constraint is active over intervals of finite duration when the splash-down is moved closer than \SI{57}{\degree}.
In particular, the heat flux peak is delayed and occurs during phase \StageFourOne{}.
Note that the duration and location of the bounded arc are very difficult to predict, but, thanks to the direct discretization method, the optimal switching structure is automatically determined and no \emph{a priori} guess is required.

\begin{figure}[h]
    \centering
    \includegraphics[width=0.7\linewidth]{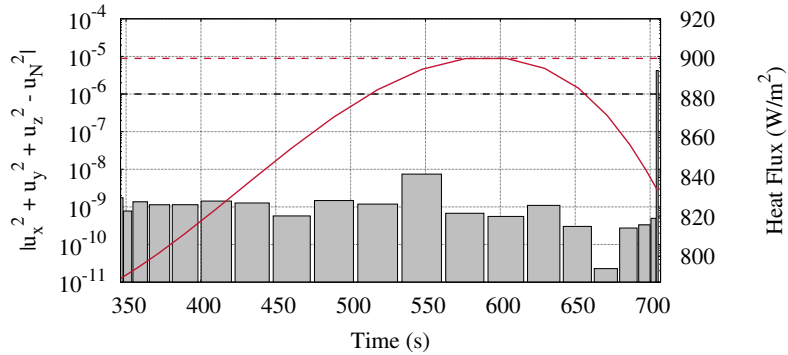}
    \caption{Relaxation error during phase \StageFourOne{} of the trajectory constrained to $\varphi_{R} = \SI{55}{\degree}$}
    \label{fig:cone_err}
\end{figure}

Fig.~\ref{fig:cone_err} reports the relaxation error during phase \StageFourOne{} of the solution corresponding to $\varphi_R = \SI{55}{\degree}$.
The error is always below the solver feasibility threshold (\smash{$10^{-6}$}), except for the final node of the phase that, due to the Radau discretization scheme, is not an optimization variable and is extrapolated \emph{a posteriori} from the approximating polynomial.
This solution is particularly interesting as  Assumption~\ref{assumption:hf} does not hold anymore in the interval $[\num{576.6}, \num{604.2}]$ \si{\s}.
Nevertheless, even though no theoretical proof can be provided, the relaxation is still lossless, as the resulting controls satisfy Eq.~\eqref{eq:thrust_direction_cone_con} with the equality sign within tolerance.

\section{Conclusion}

This paper presented a convex approach to the optimization of the ascent trajectory of a multistage launch vehicle.
The intrinsic nonconvexities of the problem have been effectively tackled via a thoughtful convexification process. 
This exploits a convenient change of variables to reduce the coupling of state and control, and it preserves some of the original problem nonlinearity by relaxing, in a lossless fashion, a control constraint.
These expedients, combined with successive linearization, are essential to set up a convex formulation of the original problem that does not suffer from numerical issues, such as high-frequency jitters in the control, and other undesired phenomena linked to the linearization that may hinder convergence.
In this respect, virtual controls and buffer zones ensure the recursive feasibility of the iterative process, 
and a simple, yet effective, method to update the reference solution based on multiple previous iterations is implemented to filter out oscillations in the search space and provide further stability to the procedure.
Moreover, it was shown that the employed $hp$ discretization scheme can, on the one hand, accurately capture the complex ascent and return dynamics, and, yet, produce a computationally efficient and sparse discrete problem.

The devised algorithm exhibits great robustness to the initialization, as it can achieve convergence even starting from a rough reference trajectory, and high overall computational efficiency, as the sequential process terminates successfully and quickly, after just a few iterations.
Thus, it represents a fast and reliable alternative to traditional optimization methods, which, in turn, often manifest high sensitivity to the supplied first guess solution or require a large computational effort to achieve convergence.
These beneficial properties make the proposed approach potentially suitable for further studies and applications to optimization-based guidance, as speed-ups can be achieved if: better initialization is provided; custom SOCP solvers are used; and the code is executed on dedicated hardware.
Naturally, specific validation tests are necessary to rigorously demonstrate the real-time applicability of the algorithm and will be the subject of future work.

In the present paper, we investigated the specific VEGA launch vehicle configuration and analyzed its performance sensitivity to the splash-down location of the third stage.
Results show that moving the return point of the spent stage can significantly change the mission profile.
As a result, the payload undergoes different, hardly predictable, thermal conditions.
Nevertheless, the convex optimization approach was proved to effectively handle both the heat flux and splash-down requirements in a systematic way, retaining system performance to acceptable levels in a wide range of scenarios.
While the numerical results are relative to this particular case study, the general approach can be easily extended to different missions and vehicle configurations.
Future work will investigate performance and flexibility of the proposed algorithm in other realistic scenarios, also including further (nonconvex) constraints to account for additional mission requirements (e.g., visibility aspects).

\section{Acknowledgment}
This work was supported by the Agreement n. 2019-4--HH.0 CUP F86C17000080005 
``Technical Assistance on Launch Vehicles and Propulsion'' 
between the Italian Space Agency and the Department of Mechanical and Aerospace Engineering of Sapienza University of Rome.
    
\bibliographystyle{AAS_publication}   
\bibliography{references_no_doi}   

\end{document}